\documentclass[11pt]{article}

\usepackage{amsmath,amsthm,amssymb,amsfonts}
\usepackage{mathtools,array,booktabs,float}
\usepackage[hidelinks,hyperfootnotes=false]{hyperref}
\newtheorem{theorem}{Theorem}[section]
\newtheorem{corollary}[theorem]{Corollary}
\newtheorem{lemma}[theorem]{Lemma}
\newtheorem{claim}[theorem]{Claim}
\newenvironment{Proof}[1][Proof]
  {\par\noindent\textbf{#1.}\ }
  {\hfill\rule{3mm}{3mm}\par\medskip}
\newcommand{\remark}{\medskip\par\noindent\textbf{Remark.}\ }

\newcommand{\Z}{\mathbb Z}
\newcommand{\R}{\mathbb R}
\DeclareMathOperator{\supp}{supp}
\newcommand{\abs}[1]{\lvert#1\rvert}
\newcommand{\JCTB}{\emph{J. Combin. Theory Ser. B}}
\newcommand{\JGT}{\emph{J. Graph Theory}}
\newcommand{\SIAMDM}{\emph{SIAM J. Discrete Math.}}

\title{Nowhere-zero $9$-flows of signed graphs}
\newcommand{\affiliationmark}[1]{%
  \hyperlink{author-affiliation-#1}{\textsuperscript{#1}}}
\author{Miaomiao Han\affiliationmark{1}, Jiaao Li\affiliationmark{2},\enspace
  Rong Luo\affiliationmark{3}, and Xiangqian Zhou\affiliationmark{4}}
\date{}
\hypersetup{pdfauthor={Miaomiao Han, Jiaao Li, Rong Luo, Xiangqian Zhou}}

\makeatletter
\newcommand{\firstpageauthorinfo}{%
  \begingroup
  \renewcommand{\thefootnote}{}%
  \long\def\@makefntext##1{\noindent##1}%
  \footnotetext{%
    \begin{minipage}[t]{\linewidth}
      \footnotesize
      \raggedright
      \hypertarget{author-affiliation-1}{\textsuperscript{1}}\enspace
      College of Mathematical Science, Tianjin Normal University, Tianjin 300387, China. \hspace*{1em}{Email:}~\texttt{mmhan2018@hotmail.com}\par\smallskip
      \hypertarget{author-affiliation-2}{\textsuperscript{2}}\enspace
      School of Mathematical Sciences and LPMC, Nankai University, Tianjin 300371, China\par
      \hspace*{1em}{Email:}~\texttt{lijiaao@nankai.edu.cn}\par\smallskip
      \hypertarget{author-affiliation-3}{\textsuperscript{3}}\enspace
      Department of Mathematics, West Virginia University, Morgantown, WV 26505, USA\par
      \hspace*{1em}{Email:}~\texttt{rluo@mail.wvu.edu}\par\smallskip
      \hypertarget{author-affiliation-4}{\textsuperscript{4}}\enspace
      Department of Mathematics and Statistics, Wright State University, Dayton, Ohio 45435\par
      \hspace*{1em}{Email:}~\texttt{xiangqian.zhou@wright.edu}
    \end{minipage}%
  }%
  \endgroup
}
\makeatother

\begin{document}
\maketitle
\firstpageauthorinfo

\begin{abstract}
In 1983, Bouchet conjectured that every flow-admissible signed graph
admits a nowhere-zero $6$-flow. DeVos et al.\ proved that every
flow-admissible signed graph admits a nowhere-zero $11$-flow. In this
paper, we show that every flow-admissible signed graph admits a
nowhere-zero $9$-flow.

\end{abstract}

\section{Introduction}
All graphs and signed graphs considered in this paper are finite and
may have multiple edges or loops. For terminology and notation not
defined here, we follow~\cite{BM2008,Diestel2010,West1996}.

Bouchet's $6$-flow conjecture~\cite{Bouchet} states that
\emph{every flow-admissible signed graph admits a nowhere-zero $6$-flow}.
Bouchet~\cite{Bouchet} proved that such graphs admit nowhere-zero $216$-flows.
Z\'yka~\cite{Zyka} improved this bound to $30$. DeVos, Li, Lu, Luo,
Zhang, and Zhang~\cite{DLLLZZ} proved that every flow-admissible signed
graph admits a nowhere-zero $11$-flow.

In this paper, we prove the following result.

\begin{theorem}\label{thm:main}
Every flow-admissible signed graph admits a nowhere-zero $9$-flow.
\end{theorem}

Several smaller bounds are known for particular classes of signed graphs.  Wu, Ye, Zang, and Zhang~\cite{Wu2014} proved  that every $8$-edge-connected flow-admissible signed graph admits a nowhere-zero $3$-flow.
Raspaud and Zhu~\cite{RZ} proved that every $4$-edge-connected
flow-admissible signed graph admits a nowhere-zero $4$-flow. It is also know that a flow-admissible signed graph admits a nowhere-zero $8$-flow if it is $3$-edge-connected as proved by DeVos, Nurse, and \v{S}\'amal~\cite{DNS}, or if it admits a nowhere-zero $\Z_4$-flow as shown by Han, Li, Lu, and Luo~\cite{HLLL},  or if its underlying graph admits a nowhere-zero integer $4$-flow as proved by Luo, M\'a\v{c}ajov\'a, \v{S}koviera, and Zhang~\cite{LMSZ}. 
Nurse~\cite{Nurse} proved Bouchet's conjecture for cyclically
$5$-edge-connected cubic signed graphs. Our Theorem \ref{thm:main} provides a new general bound for all flow-admissible signed graphs.

DeVos et al.~\cite[Theorem~1.3]{DLLLZZ} proved that every
flow-admissible signed graph admits a balanced nowhere-zero
$\Z_2\times\Z_3$-flow. A $\Z_2\times\Z_3$-flow $(f_1,f_2)$ is called
\emph{balanced} if $\supp(f_1)$ contains an even number of negative edges.
DeVos et al.~\cite{DLLLZZ} convert the $\Z_2$-flow and the $\Z_3$-flow into integer
flows to obtain a nowhere-zero $11$-flow. To prove Theorem~\ref{thm:main}, we show that such a balanced nowhere-zero
$\Z_2\times\Z_3$-flow can be converted into a nowhere-zero $9$-flow. A new ingredient in the proof is to select a half-integral solution of signed incidence
equations with prescribed endpoint coordinates and use it to construct
the required integer flow.

\section{Notation and terminology}\label{sec:prelim}

We follow the flow notation in~\cite{DLLLZZ,HLLL}. A \emph{signed graph}
$(G,\sigma)$ consists of a graph $G$ and a signature
$\sigma\colon E(G)\to\{1,-1\}$. An edge $e$ is \emph{positive} if
$\sigma(e)=1$ and \emph{negative} if $\sigma(e)=-1$. The signature is omitted when it
is clear from the context. For a subgraph $H$ of $G$, let $E_N(H)$
be its set of negative edges. A \emph{circuit} is a connected $2$-regular
subgraph, where a loop contributes two to the degree of its vertex.
 A circuit is
\emph{balanced} if it contains an even number of negative edges, and is
\emph{unbalanced} otherwise.

Every edge consists of two half-edges, each incident with one end.
The half-edges of a loop are distinct. Let $\mathcal{H}(G)$ be the set of
half-edges and let $\mathcal{H}_G(v)$ be the set of half-edges incident with
$v$. For a half-edge $a$, denote the edge containing it by $e_a$.
An orientation of $G$ is a mapping
$\tau\colon \mathcal{H}(G)\to\{1,-1\}$ such that
\[
 \tau(a)\tau(\widehat a)=-\sigma(e_a)
\]
for the two half-edges $a,\widehat a$ of each edge. A half-edge $a$
is oriented away from its end if $\tau(a)=1$, and towards its end
if $\tau(a)=-1$. Denote the resulting oriented signed graph by
$(G,\tau)$.

Fix an orientation $\tau$ of $G$ and let $B$ be its signed incidence
matrix, with
\[
 B_{v,e}=\sum_{\substack{a\in \mathcal{H}_G(v)\\e_a=e}}\tau(a).
\]
Let $A$ be an abelian group. We regard an assignment of edge values
as a vector $f\in A^{E(G)}$, with $f(e)$ denoting its coordinate at $e$.
Its \emph{boundary} under $\tau$ is the vector $Bf\in A^{V(G)}$, where
\[
 (Bf)(v)=\sum_{a\in \mathcal{H}_G(v)}\tau(a)f(e_a)
 \qquad(v\in V(G)).
\]
The vector $f$ is an \emph{$A$-flow} if $Bf=0$. Its \emph{support} is
$\supp(f)=\{e\in E(G):f(e)\ne0\}$, and $f$ is nowhere-zero if
$\supp(f)=E(G)$. An integer $k$-flow, or simply a $k$-flow, is a
$\Z$-flow $f$ with $|f(e)|<k$ for each edge $e$. A
$\Z_k$-flow is also called a modulo $k$-flow. We use $k$-NZF and
$A$-NZF for nowhere-zero $k$-flow and nowhere-zero $A$-flow,
respectively. For an integer flow $g$ and a positive integer $i$, let
\[
 E_{i}(g)=\{e\in E(G):|g(e)|=i\}.
\]

A signed graph is \emph{flow-admissible} if it admits a $k$-NZF for some
positive integer $k \geq 2$.

Reversing the directions of both half-edges of an edge and negating
the corresponding coordinate of $f$ preserves every boundary
equation. Thus a flow can be expressed as a vector with respect to
any fixed orientation.
A positive loop has a zero column in $B$, and a negative loop has
one nonzero entry, equal to $2$ or $-2$.

\section{Half-integral solutions of signed incidence equations}
\label{sec:half}

The next lemma concerns an arbitrary signed incidence matrix. For
vectors indexed by $E(G)$, inequalities are understood coordinatewise.
A \emph{vertex} of a convex set is a point that is not a nontrivial convex combination of two distinct points of that set.

\begin{lemma}\label{lem:half-integral}
Let $(G,\tau)$ be an oriented signed graph with signed incidence
matrix $B$, and let $b\in\Z^{V(G)}$. Every vertex of the polytope
\begin{equation}\label{eq:polytope}
 P_b=\{x\in\R^{E(G)}:Bx=b,\ 0\leq x\leq1\}
\end{equation}
belongs to $\{0,1/2,1\}^{E(G)}$.
\end{lemma}

\begin{Proof}
Let $x$ be a vertex of $P_b$ and let
$E_x=\{e\in E(G):0<x(e)<1\}$. An edge in $E_x$ is called a \emph{fractional edge}.

If $E_x=\varnothing$, then $x$ is integral and the conclusion holds.
Assume that $E_x\ne\varnothing$, and let $G_x$ be the subgraph with edge
set $E_x$ and vertex set consisting of the ends of these edges.

\begin{claim} \label{cl:independent}
The columns of $B$ indexed by $E_x$ are linearly independent.
Moreover, for each component $C$ of $G_x$, the columns of the submatrix
with rows indexed by $V(C)$ and columns indexed by $E(C)$ are linearly
independent.
\end{claim}
\begin{proof}
Suppose otherwise. Then there is a nonzero vector $d$ supported on
$E_x$ such that $Bd=0$. Choose
\[
 0<\epsilon<\min_{d(e)\ne0}
 \frac{\min\{x(e),1-x(e)\}}{|d(e)|}.
\]
The minimum is positive because the support of $d$ is a nonempty
subset of $E_x$. Both $x+\epsilon d$ and $x-\epsilon d$ lie in the
unit box, and $Bd=0$ gives $B(x\pm\epsilon d)=b$. They are distinct
points of $P_b$ with midpoint $x$, contrary to the choice of $x$ as
a vertex of $P_b$. This proves the first part of the claim.

The columns indexed by $E(C)$ are zero outside the rows indexed by
$V(C)$, so restricting to those rows preserves their independence.
\end{proof}

\begin{claim}\label{cl:unbalanced-circuit}
Each component of $G_x$ is an unbalanced circuit.
\end{claim}
\begin{proof}
Every coordinate outside $E_x$ is zero or one. Moving the
contributions of these coordinates to the right-hand side of $Bx=b$
gives an integral right-hand side. If $B_e$ denotes the column of $B$
indexed by $e$, then
\[
 \sum_{e\in E_x} B_e x(e)
 =b-\sum_{e\notin E_x}B_e x(e)\in\Z^{V(G)}.
\]
If a vertex of $G_x$ had degree one, its only
incident edge would be a nonloop edge $e$, and its equation would have
the form $\pm x(e)=a$ for an integer $a$. This contradicts
$0<x(e)<1$. Thus every vertex of $G_x$ has degree at least two.

Consider a connected component $C$ of $G_x$, with $n$ vertices and $m$ edges.
Its incidence columns are independent by Claim~\ref{cl:independent} and supported on its $n$ rows,
so $m\leq n$. The degree bound gives $2m\geq2n$. Consequently $m=n$ and
every vertex of $C$ has degree two. Hence $C$ is a circuit, possibly a loop or a pair of parallel edges. If $C$ consists of a loop, its column is nonzero by independence. The loop is therefore negative,
and its equation $\pm2x(e)=a$ shows that $2x(e)$ is integral.

Suppose that \(C\) has at least two edges, and let $M=B[V(C),E(C)]$.
By ordering the vertices and edges cyclically, choosing an appropriate orientation, and multiplying rows or columns by \(-1\) when necessary, we may assume that \(M\) has the form
\[
M=
\begin{pmatrix}
1      & 0      & 0      & \cdots & 0      & \varepsilon\\
-1     & 1      & 0      & \cdots & 0      & 0\\
0      & -1     & 1      & \ddots & \vdots & \vdots\\
\vdots & \ddots & \ddots & \ddots & 0      & 0\\
0      & \cdots & 0      & -1     & 1      & 0\\
0      & \cdots & \cdots & 0      & -1     & 1
\end{pmatrix},
\]
where
\[
\varepsilon=
\begin{cases}
-1, & \text{if \(C\) is balanced},\\
 1, & \text{if \(C\) is unbalanced}.
\end{cases}
\]
A direct computation yields that $\det M=1+\varepsilon$. Thus, \(\det M=0\) if \(C\) is balanced, while \(\det M= \pm 2\) if \(C\) is unbalanced. Since the columns of \(M\) are linearly independent by Claim~\ref{cl:independent}, we have \(\det M\ne0\). Therefore, \(C\) is unbalanced.
\end{proof}

\noindent
\textbf{The final step.} Let $C$ be a component of $G_x$. By
Claim~\ref{cl:unbalanced-circuit}, $C$ is an unbalanced circuit.
If $C$ consists of a loop $e$, then $2x(e)$ is integral by the
preceding argument, and hence $x(e)=1/2$.

Suppose that $C$ has at least two edges. Again, let $M=B[V(C),E(C)]$.
Moving the contributions of all edges outside $E(C)$ to the
right-hand side of the equations indexed by $V(C)$ gives $Mx_C=a$
for some $a\in\Z^{V(C)}$, where $x_C$ is the restriction of $x$ to
$E(C)$. By the proof of Claim~\ref{cl:unbalanced-circuit},
$\det M=\pm 2$. Since $M$ is an integer matrix,  $M^{-1}=\frac{1}{\det M}\operatorname{adj}(M)$ has entries in $\frac12\Z$. Thus, every coordinate of $x_C$ is a half-integer. Since each lies strictly between zero and one, each is $1/2$. Every coordinate outside $E_x$ is integral, which proves the lemma.
\end{Proof}

The proof also shows that the fractional edges form vertex-disjoint
unbalanced circuits. In particular, every bridge of $G$ has an integral
coordinate at each vertex of $P_b$.

For the flow construction, the half-integral vector must preserve selected coordinates at their lower or upper bounds. These requirements can be imposed through a face of \eqref{eq:polytope}.

\begin{corollary}\label{cor:endpoints}
Let $b\in\Z^{V(G)}$ and $z\in P_b$. Then there exists
$x\in P_b\cap\{0,\frac{1}{2},1\}^{E(G)}$ such that $z(e) = x(e)$ whenever $z(e) \in \{0,1\}$. 
\end{corollary}

\begin{Proof}
  Among all $x\in P_b$ such that $x(e) = z(e)$ for each $e$ with $z(e) \in \{0,1\}$, choose one with as many coordinates as possible equal to $0$ or $1$. This is possible since $z$ itself is such a vector. The number of coordinates of $x$ in $\{0,1\}$ is an
integer between $0$ and $|E(G)|$. By Lemma~\ref{lem:half-integral}, it is sufficient to show that $x$ is a vertex of $P_b$.

Suppose that $x=\lambda u+(1-\lambda)v$ for distinct vectors
$u,v\in P_b$ and some $0<\lambda<1$. If $x(e)$ is $0$ (resp. $1$),
then $u(e)=v(e)=x(e) =0$ (resp. $1$), since both $u(e)$ and $v(e)$ lie in $[0,1]$. Thus, $d=u-v$ is a nonzero vector with $Bd=0$, and $d(e)=0$
whenever $x(e)$ is $0$ or $1$.

Let $\delta=\min\left(
\left\{-\frac{x(e)}{d(e)}:d(e)<0\right\}
\cup
\left\{\frac{1-x(e)}{d(e)}:d(e)>0\right\}
\right)$.
Then \(\delta>0\), and for every \(0\le t\le\delta\), $0\le x+td\le1$.
Since \(Bd=0\), we also have $B(x+td)=b$, and hence,  \(x+td\in P_b\). Moreover, if \(x(e)\in\{0,1\}\), then \(d(e)=0\), so $(x+td)(e)=x(e)$.

Let $y=x+\delta d$. By the choice of \(\delta\), at least one coordinate of \(y\) that was strictly between \(0\) and \(1\) in \(x\) is now equal to \(0\) or \(1\). Thus, $y$ has at least one more coordinate equal to \(0\) or \(1\) than \(x\), contrary to our choice of $x$.

Therefore, $x$ is a vertex of $P_b$. By Lemma~\ref{lem:half-integral},
$x\in\{0,\tfrac12,1\}^{E(G)}$, as desired.
\end{Proof}

\section{Proof of the 9-flow theorem}
\label{sec:proof}

Let us recall the balanced $\Z_2\times\Z_3$-flow result of DeVos
et al.~\cite[Theorem~1.3]{DLLLZZ}.

\begin{theorem}[DeVos, Li, Lu, Luo, Zhang, and Zhang~{\cite[Theorem~1.3]{DLLLZZ}}]
\label{thm:z2z3}
Every flow-admissible signed graph admits a balanced
$\Z_2\times\Z_3$-NZF.
\end{theorem}

Since $\Z_6$ and $\Z_2\times\Z_3$ are isomorphic, we can represent a
$\Z_2\times\Z_3$-flow by a single vector over $\Z_6$. We first combine
Theorem~\ref{thm:z2z3} with the following result of Cheng et al.~\cite{CLLZ}.
\begin{theorem}[Cheng, Lu, Luo, and Zhang~{\cite[Theorem~1.5]{CLLZ}}]
\label{thm:2-to-3}
Let $(G,\tau)$ be a connected oriented signed graph with a $\Z_2$-flow $f_1$ whose
support contains an even number of negative edges. Then $G$ admits a
$3$-flow $g$ such that $E_1(g) = \supp(f_1)$. 

%$E_{g=\pm1}=\supp(f_1)$.
\end{theorem}

\begin{lemma}
\label{le:pair}
Let $(G,\tau)$ be a flow-admissible oriented signed graph with signed
incidence matrix $B$. Then there exist a vector
$f\in\{1,2,3,4,5\}^{E(G)}$ and a $3$-flow $g$ such that
$Bf\equiv0\pmod6$ and, for every edge $e$, $|g(e)|=1$ if and only if
$f(e)$ is odd.
\end{lemma}
\begin{proof}
It suffices to prove the lemma for connected $G$, since the vectors on
the connected components combine to give the desired vectors on $G$.
Let $(f_1,f_2)$ be a balanced $\Z_2\times\Z_3$-NZF with respect to $\tau$,
which exists by Theorem \ref{thm:z2z3}. Let
$f\in\{1,2,3,4,5\}^{E(G)}$ be the integer representative of
$3f_1+4f_2$ in $\Z_6^{E(G)}$. Then $Bf\equiv0\pmod6$, and $f(e)$ is
odd precisely when $e\in\supp(f_1)$.

By Theorem~\ref{thm:2-to-3}, there is a $3$-flow $g$ with respect to
$\tau$ such that $E_1(g)=\supp(f_1)=\{e\in E(G):f(e)\text{ is odd}\}$.
Thus, $f$ and $g$ satisfy the desired conditions.
\end{proof}

We call a pair $(f,g)$ satisfying Lemma~\ref{le:pair} a \emph{good
pair} for $(G,\tau)$. Theorem~\ref{thm:main} follows directly from Lemma~\ref{le:pair} and the following theorem.

\begin{theorem}
\label{thm:conversion}
Let $(G,\tau)$ be a flow-admissible oriented signed graph with signed
incidence matrix $B$, and let $(f,g)$ be a good pair for $(G,\tau)$.
Then there exists $x\in\{0,1/2,1\}^{E(G)}$ with $Bx=Bf/6$ such that
the vector
\begin{equation}\label{eq:constructed-flow}
 h=f-6x+2g,
\end{equation}
is a $9$-NZF with respect to $\tau$.
\end{theorem}

\begin{Proof}
Define
\begin{equation}\label{eq:bz}
 b=\frac{Bf}{6},\qquad z=\frac{f-g}{6}.
\end{equation}
Since $Bf\equiv0\pmod6$, the vector $b$ is integral. If $f(e)$ is odd, then
$f(e)\in\{1,3,5\}$ and $g(e)\in\{-1,1\}$. If $f(e)$ is even, then
$f(e)\in\{2,4\}$ and $g(e)\in\{-2,0,2\}$. The possible differences are
\[
\begin{array}{c|ccccc}
 f(e)&1&2&3&4&5\\ \hline
 f(e)-g(e)&0,2&0,2,4&2,4&2,4,6&4,6
\end{array}
\]
and hence, $f(e)-g(e)\in\{0,2,4,6\}$. It follows that  $z(e)\in\{0,\frac{1}{3},\frac{2}{3},1\}$ for each $e\in E(G)$.
%\begin{equation}\label{eq:z-values}
 %z(e)\in\{0,1/3,2/3,1\}\qquad(e\in E(G)).
%\end{equation}
Since $Bg=0$, we have $Bz=(Bf-Bg)/6=b$. Thus $z\in P_b$. By
Corollary~\ref{cor:endpoints}, there exists
$x\in P_b\cap\{0,1/2,1\}^{E(G)}$ such that $x(e) =z(e)$ whenever $z(e) \in \{0, 1\}$.

The vector $6x$ is integral, so $h=f-6x+2g$ is integral. Moreover,  
\begin{equation}\label{eq:conservation}
 Bh=Bf-6Bx+2Bg=6b-6b=0.
\end{equation}
It remains to show that $1\leq |h(e)|\leq8$ for every edge $e$.

For any $e\in E(G)$, since $6x(e)\in\{0,3,6\}$, we have that
\begin{equation}\label{eq:mod-three}
 h(e)=6z(e)+3g(e)-6x(e)\equiv6z(e)\pmod3.
\end{equation}
If $z(e)$ is $\frac{1}{3}$ or $\frac{2}{3}$, then $6z(e)$ is $2$ or $4$.
Neither is divisible by $3$, so $h(e)\ne0$.

If $z(e)$ is $0$ or $1$, then $x(e)=z(e)$ and hence $h(e)=3g(e)$.
In this case we claim that $g(e)$ is nonzero: otherwise $f(e)=6z(e)\in\{0,6\}$, contradicting the fact that $1\leq f(e)\leq5$.
This proves that $h$ is nowhere-zero.

For the bound, if $f(e)$ is odd, then $1\leq f(e)\leq5$ and
$g(e)\in\{-1,1\}$, and hence, $-1\leq f(e)+2g(e)\leq7$.
If $f(e)$ is even, then $2\leq f(e)\leq4$ and
$g(e)\in\{-2,0,2\}$, and hence, $-2\leq f(e)+2g(e)\leq8$.

Recall that $0\leq6x(e)\leq6$. So we have 
\[
 -8\leq h(e) = f(e)+2g(e)-6x(e)\leq8.
\]
Therefore, $1\leq\abs{h(e)}\leq8$ for every edge, as required.
\end{Proof}

\remark Table~\ref{tab:values} lists the edge values allowed by the
good-pair conditions. There are exactly two cases in which $|h(e)|=8$:
\[
 (f(e),g(e),x(e))=(2,-2,1)
 \quad\text{and}\quad
 (f(e),g(e),x(e))=(4,2,0).
\]
Consequently, the flow obtained in Theorem~\ref{thm:conversion} is an
integer $8$-flow precisely when neither case occurs on any edge.
Avoiding both cases simultaneously would require an additional
selection argument.

\begin{table}[H]
\centering
\renewcommand{\arraystretch}{1.22}
\begin{tabular}{@{}ccccc@{}}
\toprule
$f(e)$ & $g(e)$ & $z(e)$ & Permitted $x(e)$ & $h(e)$\\
\midrule
$1$ & $1$  & $0$        & $0$           & $3$\\
$1$ & $-1$ & $1/3$      & $0,1/2,1$     & $-1,-4,-7$\\
$2$ & $2$  & $0$        & $0$           & $6$\\
$2$ & $0$  & $1/3$      & $0,1/2,1$     & $2,-1,-4$\\
$2$ & $-2$ & $2/3$      & $0,1/2,1$     & $-2,-5,-8$\\
$3$ & $1$  & $1/3$      & $0,1/2,1$     & $5,2,-1$\\
$3$ & $-1$ & $2/3$      & $0,1/2,1$     & $1,-2,-5$\\
$4$ & $2$  & $1/3$      & $0,1/2,1$     & $8,5,2$\\
$4$ & $0$  & $2/3$      & $0,1/2,1$     & $4,1,-2$\\
$4$ & $-2$ & $1$        & $1$           & $-6$\\
$5$ & $1$  & $2/3$      & $0,1/2,1$     & $7,4,1$\\
$5$ & $-1$ & $1$        & $1$           & $-3$\\
\bottomrule
\end{tabular}
\caption{Edge values in the conversion of Theorem~\ref{thm:conversion}.}
\label{tab:values}
\end{table}

\section*{Declaration of AI usage}

During the preparation of this manuscript, the authors used ChatGPT 6
to explore potential approaches and refine the presentation. The
authors carefully reviewed, revised, and independently verified all
AI-assisted suggestions and arguments. The authors take full
responsibility for the mathematical correctness and content of the
manuscript.

\section*{Acknowledgments}
Miaomiao Han is partially supported by National Natural Science Foundation of China (No. 12471333). Jiaao Li is partially supported by  National Natural Science Foundation of China (No. 12571371), Natural Science Foundation of Tianjin (No. 24JCJQJC00130), and the Fundamental Research Funds for the Central Universities, Nankai University. Rong Luo is partially supported by a grant from  Simons Foundation (No. 839830).

%\clearpage

\end{document}